\documentclass[11pt,leqno]{amsart}
\usepackage[mathcal]{eucal}
\usepackage[utf8x]{inputenc}

\usepackage{amsmath,amssymb,amsthm}
\usepackage[dvipsnames]{xcolor}
\usepackage[colorlinks=true, linkcolor=MidnightBlue, citecolor=OliveGreen,pagebackref]{hyperref}
\usepackage[capitalise,nameinlink]{cleveref}
\usepackage{mathdots}
\usepackage{tikz}
\usepackage{tikz-cd}
\usetikzlibrary{arrows, automata, positioning, calc}

\newtheorem{theorem}{Theorem}[section]
\newtheorem{lemma}[theorem]{Lemma}
\newtheorem{proposition}[theorem]{Proposition}
\newtheorem{corollary}[theorem]{Corollary}

\theoremstyle{definition}
\newtheorem{definition}[theorem]{Definition}

\numberwithin{equation}{section}

\newcommand{\dottimes}[1]{\times \overset{ #1 }{\dots} \times}

\DeclareMathOperator{\Aut}{Aut}
\DeclareMathOperator{\F}{\mathbb F}
\DeclareMathOperator{\N}{\mathbb N}
\DeclareMathOperator{\Z}{\mathbb Z}
\DeclareMathOperator{\Circ}{Circ}
\DeclareMathOperator{\id}{id}
\DeclareMathOperator{\Stab}{Stab}
\DeclareMathOperator{\Sym}{Sym}
\DeclareMathOperator{\codim}{codim}

\DeclareMathOperator{\gnrc}{\diamond}
\DeclareMathOperator{\sym}{sym}
\DeclareMathOperator{\con}{con}

\makeatletter
\renewcommand*\env@matrix[1][\arraystretch]{%
  \edef\arraystretch{#1}%
  \hskip -\arraycolsep
  \let\@ifnextchar\new@ifnextchar
  \array{*\c@MaxMatrixCols c}}
\makeatother

\DeclareMathOperator{\qund}{\quad\text{ and }\quad}

\newcommand{\GGS}{{GGS}}

\title[Derived series of \GGS-groups]{The derived series of \GGS-groups}

\author[J. M. Petschick]{Jan Moritz Petschick}
\address{Jan Moritz Petschick: Mathematisches
  Institut, Heinrich-Heine-Universit\"at, 40225 D\"usseldorf, Germany}
\email{jan.petschick@hhu.de}
\thanks{The research was funded by the Deutsche Forschungsgemeinschaft (DFG, German Research Foundation) — 380258175}


\begin{document}

\begin{abstract}
	Given a \GGS-group $G$ with non-constant defining tuple over a prime-regular rooted tree, we calculate the indices $|G:G^{(n)}|$ and describe the structure of the higher derived subgroups~$G^{(n)}$ for all $n \in \N$. We find that the values $|G:G^{(n)}|$ depend only mildly on the structure of the defining tuple.
\end{abstract}

\maketitle


\section{Introduction}

The class of groups of automorphisms of regular rooted trees provides many examples with interesting asymptotic and structural properties. One particularly well-studied case is the family of Grigorchuk--Gupta--Sidki-groups (usually abbreviated as `\GGS\nobreakdash-groups'). It contains at least one group of intermediate growth \cite{FG85} and many finitely generated infinite periodic groups, cf.\ \cite{GS83}. \GGS-groups are groups of automorphisms of the $p$-regular rooted tree, for an odd prime $p$, and generalise the Gupta--Sidki~$p$-groups. They are easily defined by a non-zero element $\mathbf{e}$ of $\F_p^{p-1}$ as `input data', and many of their properties can be read off this element $\mathbf{e}$: One can determine whether the corresponding \GGS-group is periodic or contains elements of infinite order, if the group is just-infinite, cf.\ \cite{Vov00}, whether it is a branch group, cf.\ \cite{FZ13}, or if it has the congruence subgroup property, cf.\ \cite{FGU17}, one may compute its Hausdorff dimension, cf.\ \cite{FZ13}, or decide if two \GGS-groups are isomorphic, cf.\ \cite{Pet19}, just by considering the defining tuples. Most of these results require subtle insights into the structure of a general \GGS-group, and some involve heavy computation. Many of the results extend to larger classes of groups, cf.\ for example \cite{AKT16, BS01, Pet21}, but have been established first for \GGS-groups, making the class of \GGS-groups a fertile soil for establishing new techniques.

However, other questions remain open. In contrast to features related to the action on the tree, many purely algebraic properties of \GGS-groups are not well-understood. In this work, we describe the derived series~$(G^{(n)})_{n \in \N}$ of all \GGS-groups, excluding those that arising from constant tuples, i.e.\ elements of the form $(\lambda, \lambda, \dots, \lambda) \in \F_p^{p-1}$ for some~$\lambda \in \F_p^\times$. A description of the derived series has previously been obtained for the special case of the Gupta--Sidki~$3$-group $\ddot \Gamma$ by Vieira in~\cite{Vie98}, along with some results concerning the lower central series of $\ddot \Gamma$. The proof, however, does not carry over to general \GGS-groups.

We now state our main result.
\begin{theorem}\label{DerGGS:thm:main}
	Let $p$ be an odd prime and let $G$ be a \GGS-group acting on a $p$-regular tree with non-constant defining tuple $\mathbf{e} \in \F_p^{p-1}$. Denote by $\mathbf{e'}$ the tuple of differences between the entries of $\mathbf{e}$, and by $\mathbf{e''}$ the tuple of differences of $\mathbf{e'}$. Then
	\[
		\log_p |G:G^{(n)}| = \begin{cases}
			p^{n-2}(p + \con(\mathbf{e'}) + \sym(\mathbf{e''})) - \frac{p^{n-1}-1}{p-1}\sym(\mathbf{e}) + 1 &\text{ if }n \geq 2,\\
			2 &\text{ for }n = 1,
		\end{cases}
	\]
	where
	\[\begin{array}{lcr}
		\sym(\mathbf d) = \begin{cases}
			1 &\text{ if } \mathbf d \text{ is symmetric,}\\
			0 &\text{ otherwise,}
		\end{cases}
		&\text{ and }&
		\con(\mathbf d) = \begin{cases}
			1 &\text{ if } \mathbf d \text{ is constant,}\\
			0 &\text{ otherwise.}
		\end{cases}
	\end{array}\]
\end{theorem}
A tuple is called \emph{symmetric} if its $i$\textsuperscript{th} entry is equal to its $i$\textsuperscript{th}-to-last entry. It is no surprise that the vector $\mathbf{e'} \in \F_p^{p-2}$ of differences between neighbouring entries in $\mathbf{e}$ is associated to the determination of the structure of the derived subgroups, since it describes the sections of the commutator $[b,a]$ of the two generators of a \GGS-group. Interestingly, the indices of the derived subgroups do not depend on the higher iterates of the differences.

It is worthwhile to compare our result with the main result of \cite{FZ13}, where the indices of the congruence subgroups, i.e.\ the pointwise stabilisers $\Stab_G(n)$ of all elements of a given distance $n \in \N$ from the root of the tree, are computed to be
\[
	\log_p |G:\Stab_G(n)| = \begin{cases}
		tp^{n-2} + \frac{p^{n-2}-1}{p-1}\sym(\mathbf{e}) + 1 &\text{ if }n \geq 2,\\
		1 &\text{ for }n = 1,
	\end{cases}
\]
where $t$ is the rank of a certain matrix associated to $\mathbf{e}$, which takes values in $\{ 2, \dots p\}$. In particular, the number of configurations of the indices $|G:\Stab_G(n)|$ grows linearly with~$p$. In comparison, the indices of the derived series are more uniform and depend only on three (interconnected) binary invariants of $\mathbf{e}$; hence the indices of the derived subgroups of any \GGS-groups (aside from the dependency on the prime $p$ itself), fall in precisely four distinct classes.

From a group-theoretic standpoint, it is an inherently interesting problem to determine the derived series of a given intriguing group. This is especially true since \GGS-groups are hypoabelian, i.e.\ the intersection of all members of the derived series is trivial; hence every element of $G$ appears as a non-trivial element in some quotient $G^{(n)}/G^{(n+1)}$. Furthermore, the derived series fulfils the analogue of the congruence subgroup property: all finite index subgroups contain some derived subgroup. This is an immediate consequence of the congruence subgroup property of \GGS-groups, that was established in \cite{FGU17}, and the fact that the $n$\textsuperscript{th} derived subgroup $G^{(n)}$ is contained in the $n$\textsuperscript{th} level stabiliser $\Stab_G(n)$. The analogy to the congruence subgroups goes further. We prove the following theorem.
\begin{theorem}\label{DerGGS:thm:higher derived subgps}
	Let $G$ be a \GGS-group with defining tuple $\mathbf{e}$ and let $n \in \N_{\geq 3}$. Then
	\[
		\psi(G^{(n)}) = G^{(n-1)} \dottimes{p} G^{(n-1)}.
	\]
	If $\con(\mathbf{e'}) + \sym(\mathbf{e''}) - \sym(\mathbf{e}) = 0$, the same holds for $n = 2$. 
\end{theorem}
Note that the equivalent statement is true for the congruence subgroups of \GGS-groups. To make the connexion between congruence and derived subgroups more transparent, we introduce the \emph{series of iterated local laws}. It is a descending series $(\mathrm L_n(G))_{n \in \N}$ of normal subgroups of a group $G$ acting on a rooted tree, such that $\mathrm L_n(G) \leq \Stab_G(n)$ for all $n \in \N$, and is formed by the elements that have to stabilise vertices of a certain distance by virtue of fulfilling certain algebraic equations in $G$. In this sense, it is the `algebraic analogue' of the sequence of layer stabilisers. See \cref{DerGGS:def:iterated local laws} for a precise definition.

As a corollary to \cref{DerGGS:thm:main}, we prove that the series of iterated local laws and the derived series coincide for all \GGS-groups defined by a non-constant vector $\mathbf{e}$, thus explaining, at least heuristically, the similarities mentioned above; see \cref{DerGGS:cor:Derived and iterated local}.

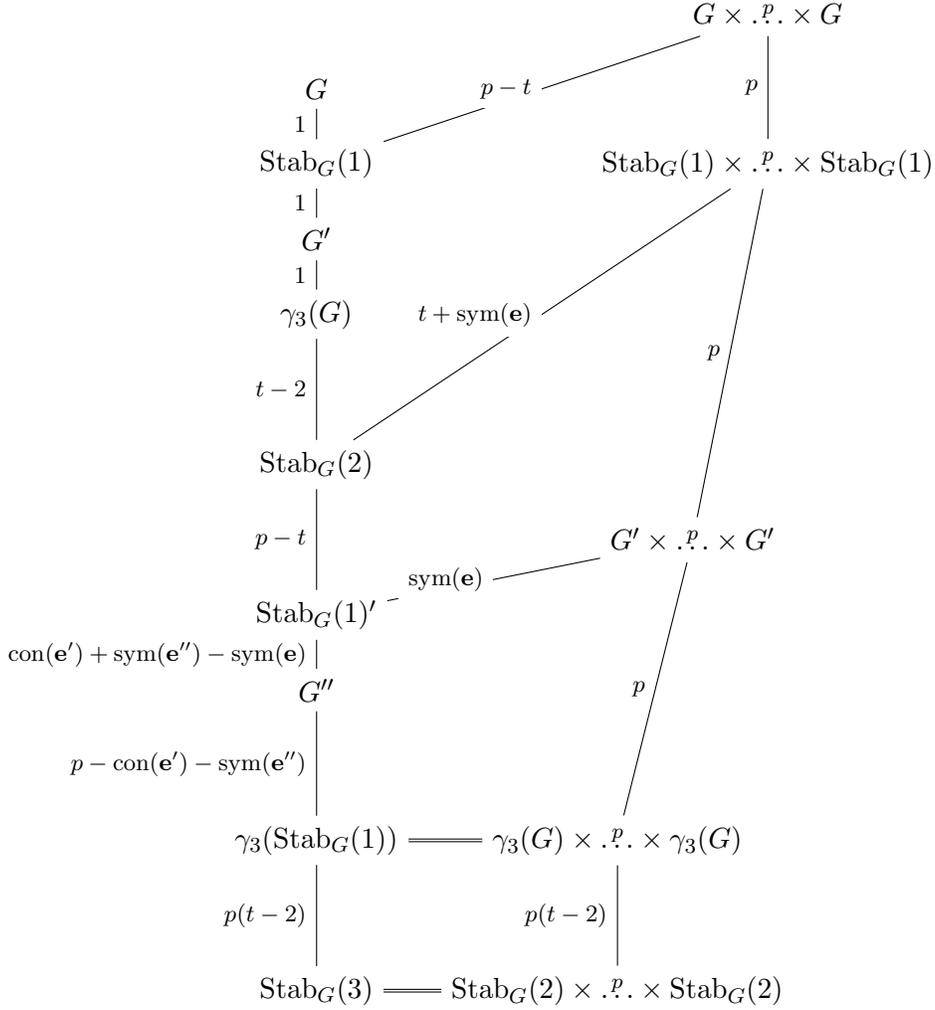
\begin{figure}
	\begin{tikzpicture}[scale=1]
		\node (G)		at	(0,12)	{$G$};
		\node (S1)		at	(0,11)	{$\Stab_G(1)$};
		\node (der)		at	(0,10)	{$G'$};
		\node (gam)		at	(0,9)	{$\gamma_3(G)$};
		\node (S2)		at	(0,7)	{$\Stab_G(2)$};
		\node (derS1)	at	(0,5)	{$\Stab_G(1)'$};
		\node (der2)	at	(0,4)	{$G''$};
		\node (gamS1)	at	(0,2)	{$\gamma_3(\Stab_G(1))$};
		\node (S3)		at	(0,0)	{$\Stab_G(3)$};
	
		\node (Gp)		at	(6,13)	{$G \dottimes{p} G$};
		\node (S1p)		at	(6,11)	{$\Stab_G(1) \dottimes{p} \Stab_G(1)$};
		\node (derp)	at	(5,6)	{$G' \dottimes{p} G'$};
		\node (gamp)	at	(4,2)	{$\gamma_3(G) \dottimes{p} \gamma_3(G)$};
		\node (S2p)		at	(4,0)	{$\Stab_G(2) \dottimes{p} \Stab_G(2)$};

		\path
			(G)			edge	node[left]	{\footnotesize{$1$}}	(S1)
			(S1)		edge	node[left]	{\footnotesize$1$}	(der)
			(der)		edge	node[left]	{\footnotesize$1$}	(gam)
	        (gam)		edge    node[left]	{\footnotesize$t-2$}	(S2)	
	        (S2)		edge    node[left]	{\footnotesize$p-t$}	(derS1)
	        (derS1)	    edge    node[left]	{\footnotesize$\con(\mathbf{e'}) + \sym(\mathbf{e''}) - \sym(\mathbf{e})$}	(der2)
	        (der2)	    edge    node[left]	{\footnotesize$p - \con(\mathbf{e'}) - \sym(\mathbf{e''})$}	(gamS1)
	        (gamS1)	    edge	node[left]	{\footnotesize$p(t-2)$}	(S3)
			(Gp)		edge	node[left]	{\footnotesize$p$}	(S1p)	
	        (S1p)		edge	node[left]	{\footnotesize$p$}	(derp)
	        (derp)	    edge	node[left]	{\footnotesize$p$}    (gamp)
			(gamp)	    edge	node[left]	{\footnotesize$p(t-2)$}    (S2p)
			(S1)		edge	node[left, fill=white]	{\footnotesize$p-t$}	(Gp)
			(S2)		edge	node[left, fill=white]	{\footnotesize$t + \sym(\mathbf{e})$}	(S1p)
			(derS1)		edge	node[left, fill=white]	{\footnotesize$\sym(\mathbf{e})$}	(derp)
			(gamS1)		edge[double]			(gamp)
			(S3)		edge[double]			(S2p);
	\end{tikzpicture}
	\caption{Part of the top of the subgroup lattice of a \GGS-group, with some supergroups added. Passage from the left to the right side signifies the application of $\psi$. All indices are logarithmic.}
	\label{DerGGS:fig:lattice}
\end{figure}

The paper is organised in the following way. After establishing our notation, we prove some structural results on \GGS-groups. Then we prove \cref{DerGGS:prop:G'' index}, in which we compute the index of the second derived subgroup in the full group. This is the main technical step. Afterwards, we proceed to derive our other results. Being aware of the multitude of subgroups appearing, we point the reader to \cref{DerGGS:fig:lattice}, which depicts the relevant portion of the top of the subgroup lattice of a \GGS-groups.

\subsection*{Acknowledgements} This is part of the author's Ph.D.~thesis, written under the supervision of Benjamin~Klopsch at the Heinrich-Heine-Universität~Düsseldorf. The author thanks Karthika~Rajeev for helpful discussions and bringing his attention to the problem, and Gustavo~Fern\'andez-Alcober, Mikel~Garciarena~Perez, Margherita~Piccolo and Djurre~Tijsma for their respective comments on preliminary versions of this paper.

\section{On Grigorchuk--Gupta--Sidki-groups}

We begin with some generalities. We fix an odd prime $p$. Given a group $G$ and two elements $g, h$, we use the following conventions for conjugation and the commutator
\[
	g^h = h^{-1}gh \qund [g,h] = g^{-1}h^{-1} g h = (h^{-1})^g h.
\]

\subsection{Groups of automorphisms of regular rooted trees}
Write~$X$ for the set $\{0, \dots, p-1\}$. We will sometimes identify $X$ with the set underlying the field~$\F_p$. We write $X^*$ for the Cayley graph of the free monoid on $X$, which is a rooted $p$-regular tree, i.e.\ a loop-free graph in which all but one vertex have valency $p+1$, and the remaining vertex $\varnothing$, called the \emph{root} of the tree, has valency $p$. The vertices of $X^*$ are the set of finite sequences in $X$. We write $X^n$ for the set of all vertices of a given length $n \in \N$, and call this set the \emph{$n$\textsuperscript{th} level of $X^*$}. The root $\varnothing$ has length $0$.

Any (graph) automorphism $g \in \Aut(X^*)$ necessarily fixes $\varnothing$, since it has fewer neighbours than every other vertex, and must consequently leave the levels $X^n$ invariant for all $n \in \N$. We write $\Stab(n)$ for the stabiliser of $X^n$, and $\Stab_G(n)$ for its intersection with some subgroup $G \leq \Aut(X^*)$. Let $u$ and $v$ be vertices of $X^*$. We write $u^g$ for the image of $u$ under $g$. Since levels are invariant under \(g\), the equation
\[
	(uv)^g = u^g v^{g|_u}
\]
uniquely defines a map $|_u: \Aut(X^*) \to \Aut(X^*)$, called the \emph{section map at \(u\)}. The image is consequently called \emph{the section of $g$ at $u$}. Using these images, any tree automorphism $g$ can be decomposed into the sections prescribing the action at the subtrees of the first level, and the \emph{action of $g$ at the root} $g|^\varnothing \in \Sym(X)$, which is defined as the action of $g$ on the first level $X = X^1$. In particular, the map
\begin{align*}
	\psi: \Stab(1) &\to \Aut(X^*) \dottimes{p} \Aut(X^*)\\
	g &\mapsto (\gnrc: g|_{\gnrc})
\end{align*}
is a group isomorphism. Here we adopt the convention that the expression
\[
	(i_0: a_0,\, \dots,\, i_k: a_k,\, \gnrc: a_{\gnrc})
\]
denotes the tuple indexed by $X$, with the object $a_0$ at position $i_0$, the object $a_1$ at position $i_1$ and so forth, and the object $a_{\gnrc}$ (maybe varying in $\gnrc$) at all other positions $\gnrc \in X\smallsetminus\{i_m \mid m = 0, \dots, k\}$. The symbol $\gnrc$ will be reserved for this use. An automorphism with at most one non-trivial section is called \emph{rooted}. Rooted automorphisms must necessarily permute the subtrees $\{xX^*\mid x \in X \}$ of the first level and can be identified with permutations of $X$.

We record some equations for sections. Let $u$ and $v$ be vertices of $X^*$ and $g$ and $h$ be any automorphisms, then
\begin{align*}
	(g|_u)|_v = g|_{uv},\quad
	(gh)|_u = g|_uh|_{u^g},\quad
	g^{-1}|_u = (g|_{u^{g^{-1}}})^{-1}.
\end{align*}

A subgroup $G \leq \Aut(X^*)$ is called \emph{self-similar}, if for all vertices $u \in X^*$, the image of the section map $|_u: g \mapsto g|_u$ is contained in $G$. A self-similar group $G$ is called \emph{contracting}, if there exists a finite set $\mathcal N \subseteq G$, such that for all $g \in G$ there exists some $n \in \N$ such that for all $m \geq n$ and all $v \in X^m$ the section $g|_v$ is an element of the finite set $\mathcal N$. For a contracting group, there is a unique minimal set $\mathcal N$ with this property, which is called the \emph{nucleus of $G$}. A group $G \leq \Aut(X^*)$ is called \emph{fractal}, if for every $g \in G$ and every $x \in X$ there is an element $\widehat{g}$ stabilising the vertex $x$ such that $\widehat{g}|_x = g$. A group $G \leq \Aut(X^*)$ is called \emph{spherically transitive} if it acts transitively on every level $X^n$.

A self-similar group $G \leq \Aut(X^*)$ is called a \emph{regular branch group}, if it is spherically transitive, and if there is a finite index subgroup $K \leq G$ such that
\[
	K \dottimes{p} K \leq \psi(K).
\]
A standard technique for establishing that a group is regular branch is given by the following lemma, cf.\ \cite[Proposition~2.18]{FZ13}.
\begin{proposition}\label{DerGGS:prop:normal generation}
	Let $G \leq \Aut(X^*)$ be a spherically transitive fractal group, $H \leq G$ a subgroup and let $S \subseteq G$ be a subset. If $\{ (0: s, \; \gnrc: \id) \mid s \in S \}$ is contained in $\psi(H)$, then
	\[
		\langle S \rangle^G \dottimes{|X|} \langle S \rangle^G \leq \psi(H^G).
	\]
\end{proposition}
We now come to the definition of the series of iterated local laws for a spherically transitive group $G \leq \Aut(X^*)$. Consider the set $\{ (g|_u)|^{\varnothing} \mid g \in G, u \in X^* \} \subseteq \Sym(X)$. In view of the equations for sections above, it is easily seen that this set forms a subgroup $P(G) \leq \Sym(X)$ of the symmetric group on $X$. Given a subgroup $H \leq \Sym(X)$, we may define a corresponding subgroup of $\Aut(X^*)$ by
\[
	\Lambda(H) = \{ g \in \Aut(X^*) \mid (g|_u)|^{\varnothing} \in H \text{ for all }u \in X^* \}.
\]
If $H$ is a $p$-group, it is necessarily cyclic, and $\Lambda(H)$ is a Sylow pro-$p$-subgroup of $\Aut(X^*)$.

Let $R = R(H)$ be the collection of all group laws of $H$, i.e.\ elements of the free group $F_\infty$ on infinitely many generators that evaluate to the trivial element for all assignments of the generators to elements of $H$. Given any group $K$, we may consider the subgroup $\mathrm L_R(K)$ generated by all verbal subgroups of $K$ corresponding to elements in $R$, i.e.\ by all images of $R$ under any assignment of the generators of $F_\infty$ to elements of $K$.

Returning to the subgroup $P(G)$ defined by a group of tree automorphisms, we see that $\mathrm L_{R(P(G))}(G)$ is contained in the first level stabiliser, since $\mathrm L_{R(P(G))}(G|^\varnothing)$ is trivial by construction; note that the map $|^\varnothing: G \to P(G) \leq \Sym(X)$ is a (not necessarily surjective) homomorphism. Consider the iterates $R_n$ of $R$, that are recursively defined by
\[
	R_n = \left\{ s(r_1, \dots, r_n) \;\middle|\; \begin{array}{l}
			s \in R_{n-1}, r_i \in R \text{ for } i \in \{1, \dots, n\},\\
			s \text{ an element involving $n$ generators, and }\\
			r_1, \dots, r_n\text{ share no generator of }F_\infty
		\end{array}\right\},
\]
for $n > 1$ and $R_1 = R$, where by the expression $s(r_1, \dots, r_n)$ we mean the element of $F_\infty$ obtained by replacing the $n$ generators occurring in $s$ by the elements $r_1$ to $r_n$, these are a generalised form of the \emph{S-type iterated identities} defined in~\cite{Ers15}. By the same argument as above,
\[
	\mathrm L_{R_n(P(G))}(G) \leq \Stab_G(n).
\]
\begin{definition}\label{DerGGS:def:iterated local laws}
	Given a spherically transitive group of $G \leq \Aut(X^*)$, we write $\mathrm L_n(G)$ for $\mathrm L_{R_n(P(G))}(G)$. We call $(\mathrm L_n(G))_{n \in \N}$ the \emph{series of iterated local laws of $G$}.
\end{definition}
We consider this series in the case of \GGS-groups. It will be apparent from the definition that every \GGS-group $G$ acts locally by permutations from a cyclic group of order~$p$, i.e. it fulfils
\[
	P(G) = \langle \,(0 \, 1 \, \dots \, p-1)\, \rangle,
\]
The laws of such a group are generated by the commutators and~$p$\textsuperscript{th} powers of generators in $F_\infty$. Consequently, the group $\mathrm L_1(G)$ is the subgroup generated by $G'$ and the $p$\textsuperscript{th} powers in $G$. We shall prove that $\mathrm L_n(G)$ is in fact equal to $G^{(n)}$.

\subsection{\GGS-groups and their defining tuples}\label{sec:def tuples}

Let
\(
	\mathbf{e} = (\mathbf e_1, \dots, \mathbf e_{p-1}) \in \F_p^{p-1}
\)
be a non-\emph{constant} tuple, i.e.\ such that there are at least two different entries. We call the group $G_{\mathbf e}$ generated by the rooted automorphism $a = (0 \; 1 \; \dots \; p-1)$ and the automorphism defined by
\[
	b = \psi^{-1}(0: b,\; \gnrc: a^{e_{\gnrc}})
\]
the \emph{\GGS-group defined by $\mathbf{e}$}, and we call $\mathbf{e}$ the \emph{defining tuple of $G_\mathbf{e}$}.

Note that we exclude all constant tuples (in particular the zero tuple). The groups defined by constant non-zero tuples in the same fashion as above are usually also referred to as \GGS-groups. Furthermore, groups defined by the same construction using elements of $\Z/m\Z$ whose entries are set-wise coprime are sometimes also referred to as \GGS-groups. In general, the structure of these groups is much less understood than in the case we consider here. Even for prime powers $m = p^n$, the situation is much more involved, see for example \cite{DFG22}, where the branching structures for these groups have been evaluated.

We consider the derived subgroups of a \GGS-group $G = G_\mathbf{e}$. Since $G$ is two-generated, the first derived subgroup is normally generated by the commutator $c = [b, a]$, whose action on the tree is given by
\begin{align*}
	\psi([b, a]) &= \psi(b^{-1}) \psi(b^a) \\
	&= (0: b^{-1},\; \gnrc: a^{-\mathbf{e}_{\gnrc}}) (1: b,\; \gnrc: a^{\mathbf{e}_{\gnrc-1}}) \\
	&= (0: b^{-1}a^{\mathbf{e}_{p-1}},\; 1: a^{-\mathbf{e}_{1}}b,\; \gnrc: a^{\mathbf{e}_{\gnrc-1} - \mathbf{e}_{\gnrc}}).
\end{align*}
This signifies the importance of the \emph{first difference tuple of $\mathbf{e}$}, which we define as
\begin{align*}
	\mathbf{e'} = (\mathbf{e'}_2, \mathbf{e'}_3, \dots, \mathbf{e'}_{p-2}, \mathbf{e'}_{p-1}) \in \F_p^{p-2},
\end{align*}
where $\mathbf{e'}_i = \mathbf{e}_{i-1}-\mathbf e_i$ for all $i \in \{2, \dots, p-1\}$. We shall see that the index of the second derived subgroup in $G$ depends furthermore on the \emph{second difference tuple of $\mathbf{e}$}, given by
\[
	\mathbf{e''} = (e''_3, e''_4, \dots, e''_{p-2}, e''_{p-1}) \in \F_p^{p-3},
\]
where $\mathbf{e''}_i = \mathbf{e}_{i-2} - 2\mathbf{e}_{i-1} + \mathbf{e}_i = \mathbf{e'}_{i-1} - \mathbf{e'}_i$ for all $i \in \{3, \dots, p-1\}$. In case $p = 3$ the tuple $\mathbf{e''}$ is the empty tuple. Clearly, we have described the beginning of an iterative procedure, but
surprisingly, the indices of the higher derived subgroup do not depend on `higher' difference tuples.

In the following, we consider the elements of the vector space $\F_p^p$ and more generally of direct products of groups $G \dottimes{p} G$ as indexed by the set $X$. The choice of indexing for the defining tuple and its differences we have made above is for the following reason. Let $\log_a: \langle a \rangle \to \F_p$ be the map assigning the power $(i \bmod p)$ to any $a^i$, and $\theta: G \to \F_p^p$ the map $g \mapsto (\log_a(g|_0)^\varnothing, \dots, \log_a(g|_{p-1})^{\varnothing})$ assigning to $g$ its local actions under the first layer vertices. Then, by definition,
\[
	\theta(b) = (0, \mathbf e_1, \dots, \mathbf e_{p-1}).
\]
Thus we think of the defining tuple as an `incomplete element' of $\F_p^p$, and the element above as its full counterpart; similarly we think of $\mathbf{e'}$ and $\mathbf{e''}$ as the `tails' of regularly formed elements of
\[
	\theta(c) = (\mathbf e_{p-1}, -\mathbf e_1, \mathbf{e'}_2, \mathbf{e'}_3 \dots, \mathbf{e'}_{p-1})
\]
and
\[
	\theta([c, a]) = (\mathbf e_{p-2}-2\mathbf e_{p-1}, \mathbf e_1 + \mathbf e_{p-1}, -2\mathbf e_1+\mathbf e_2, \mathbf{e''}_3, \mathbf{e''}_4, \dots, \mathbf{e''}_{p_1})
\]
respectively.

We call $\mathbf{e}$, resp.\ $\mathbf{e''}$, \emph{symmetric} if and only if
\begin{align*}
	\mathbf e_i &= \mathbf e_{p-i} &\text{ for all } i \in \{1, \dots, p-1\}, \text{ resp. }\\
	\mathbf{e''}_{i+1} &= \mathbf{e''}_{p+1-i} &\text{ for all } i \in \{2, \dots, p-2\}.
\end{align*}
These are clearly linear conditions. We define two linear subspaces $S = \ker(M)$ and $\ddot{S} = \ker(\ddot{M})$ of $\F_p^{p-1}$ as the kernels of the two linear maps given by the matrices
\[
	M = \begin{pmatrix}
		1 & 0 & 0 & \dots & 0 & 0 & -1\\
		& 1 & 0 & \dots & 0 & -1 &\\
		& & \ddots & & \iddots &\\
		& & & 1 \;\; -1 & &
	\end{pmatrix}^\mathrm{t},
\]
in $\operatorname{Mat}(p-1, (p-1)/2; \F_p)$ and
\[
	\ddot{M} = \begin{pmatrix}
		1 & -2 &  1 &      0 &  0 & \dots & \dots &  0 &       0 & -1 &  2 & -1\\
		  &  1 & -2 &      1 &  0 & \dots & \dots &  0 &      -1 &  2 & -1 &\\
		  &    &    & \ddots &    &       &       &    & \iddots &    &    &\\
		  &    &    &      1 & -2 &     1 & -1    &  2 &      1  &    &\\
		  &    &    &        &  1 &    -3 &  3    & -1 &         &    &
	\end{pmatrix}^\mathrm{t}
\]
in $\operatorname{Mat}(p-1, (p-3)/2 ;\F_p)$. Clearly $\mathbf{e}$ is symmetric if and only if $\mathbf{e}$ is an element of $S$, and we have $\mathbf{e} \in \ddot{S}$ if and only if $\mathbf{e''}$ is symmetric. Using this description, the following lemma becomes apparent.

\begin{lemma}\label{DerGGS:lem:e sym implies snd der sym}
	Let $\mathbf{e}$ be symmetric. Then $\mathbf{e''}$ is symmetric.
\end{lemma}

\begin{proof}
	It can be easily seen that the subspace generated by the columns (displayed above as rows) of $\ddot{M}$ is contained in the subspace generated by the columns of $M$, i.e.\ that $S \subseteq \ddot{S}$. Thus all symmetric $\mathbf{e}$ yield symmetric $\mathbf{e''}$.
\end{proof}

Of course, if a vector is constant, all its difference tuples are trivial, hence in particular constant and symmetric.

We can also see that the containment of \cref{DerGGS:lem:e sym implies snd der sym} is proper and that $\codim_{\ddot{S}} S = 1$. For further computations, we simplify the basis given by the columns of $\ddot{M}$ using Gau{\ss}-Jordan elimination, and obtain
\[
	\ddot{N} = \begin{pmatrix}
			1 & 0 &  \dots & 0 &      2 &     -2 &  0 &   \dots & 0  & -1\\
			  & 1 &  \dots & 0 &      4 &     -4 &  0 &   \dots & -1 &\\
			  &   & \ddots &   & \vdots & \vdots &    & \iddots &    &\\
			  &   &        & 1 &     -3 &      3 & -1 &         &    &
	\end{pmatrix}^\mathrm{t}.
\]

\begin{lemma}\label{DerGGS:lem:linear eq a consequence of snd der sym}
	Let $\mathbf e \in \F_p^{p-1}$. If the second difference tuple $\mathbf{e''}$ is symmetric, then
	\[
		2(\mathbf{e}_{p-1} - \mathbf{e}_{1}) + (\mathbf{e}_{2} - \mathbf{e}_{p-2}) = 0.
	\]
\end{lemma}

\begin{proof}
	Since $\mathbf{e''}$ is symmetric, the vector $\mathbf{e}$ is contained in $\ddot{S}$. But the given linear equation is a linear combination of the first two columns of $\ddot{N}$, from whence the equality follows.
\end{proof}

For convenient use in formulas, we define the shorthand notation
\begin{align*}
	\con(\mathbf{d}) = \begin{cases}
		1 &\text{ if }\mathbf{d} \text{ is constant},\\
		1 &\text{ otherwise.}
	\end{cases} \qund \sym(\mathbf{d}) = \begin{cases}
		1 &\text{ if }\mathbf{d} \text{ is symmetric},\\
		1 &\text{ otherwise.}
	\end{cases}
\end{align*}

Aside from the first and second difference tuples and the defining tuple itself, all cyclic shifts of $\mathbf{e}$ (under the action of $a$) influence the structure of the \GGS-group $G_\mathbf{e}$. Fern\'andez-Alcober and Zugadi-Reizabal \cite{FZ13} demonstrated that the index of the level stabilisers in a \GGS-group depends only on whether the defining tuple is symmetric, and the rank of the \emph{circulant matrix $\Circ(0 \; \mathbf{e}) \in \operatorname{Mat}(p,p; \F_p)$} associated to the `full version' $\theta(b)$ of the defining tuple $(0 \; \mathbf{e})$, that is the matrix whose rows are the cyclic shifts of the vector $\theta(b)$.

For the computation of the indices of the level stabilisers and the derived subgroups, we use the König-Rados theorem, cf.\ e.g.\ \cite[\textsection~134]{Kro01} or \cite{Ing56}, which solves the problem of determining the rank of a circulant matrix over a prime field.

\begin{theorem}[König-Rados]\label{DerGGS:thm:KoenigRados}
	Let $p$ be a prime and $\mathbf{d} \in \F_p^n$ a vector. Then
	\[
		\operatorname{rk} \Circ(\mathbf{d}) = n - m,
	\]
	where $m$ is the multiplicity of $1$ as a root of the polynomial $E_\mathbf{d} = \sum_{i = 0}^{n-1} {d}_i X^{i}$. In particular, $\operatorname{rk} \Circ(\mathbf{d}) = n$ if and only if $\sum_{i = 0}^{n-1} {d}_i \neq 0$.
\end{theorem}

For our purposes, we make a more general definition. Let $V$ be a finite-dimensional vector space over a finite field, let $\mathcal{B}$ be an (ordered) basis for $V$, and let $C$ be the linear map that cyclically permutes the basis elements. Given a subset $M \subseteq V$, we denote by $\Circ(M)$ the smallest $C$-invariant subspace containing $M$, the \emph{circulant space of~$M$}. The notational conflict with the definition of the circulant matrix given above is negligible; the circulant space is just the row space of the circulant matrix. We will often make no a distinction between the two.

There are not many $C$-invariant subspaces. This is no surprise, since $C$ defines the regular representation of a group of order $p$, which is the sum of $p$ one-dimensional irreducible sub-representations. In fact, there is a unique (full) flag of $C$-invariant subspaces in $V$, which we record in the following proposition.

\begin{proposition}\label{DerGGS:prop:flag of circulant spaces}
	Let $V$ be an $\F_p$-vector space of dimension $n \in \N$, with basis $\mathcal{B} = \{ b_0, \dots, n-1\}$. Then the set of circulant spaces of $V$ has cardinality $n + 1$ and forms a full flag of $V$. In particular, for any $M \subseteq V$,
	\[
		\Circ(M) = \bigcup_{m \in M} \Circ(m).
	\]
	Furthermore, the dimension of $\Circ(m)$ for some $m = \sum_{i = 0}^{n-1} m_i b_i = m \in V$ is $i$ if and only if $i$ is minimal such that
	\[
		\sum_{j = n-i}^{n-1} \binom j {n-i} m_i = 0.
	\]
\end{proposition}

\begin{proof}
	We prove that, for every $i \in \{ 0, \dots, n \}$, the set
	\[
		\Circ_i(V) = \{ \mathbf{d} \in V \mid \operatorname{rk} \Circ(\mathbf{d}) \leq i \}
	\]
	is an $i$-dimensional $C$-invariant subspace of $V$. If this is true, the circulant space associated to any $\mathbf{d} \in \Circ_i(V) \smallsetminus \Circ_{i-1}(V)$ is in fact equal to $\Circ_i(V)$; since it is the minimal invariant subspace containing $\mathbf{d}$, it is contained in $\Circ_i(V)$, and by the choice of~$\mathbf{d}$ it has the same dimension as $\Circ_i(V)$. Now for any subset $M \subseteq V$, the circulant space is equal to the smallest invariant subspace containing all $\Circ(\mathbf{m})$ for $\mathbf{m} \in M$. It is easy to see that these spaces are linearly ordered, hence $\Circ(M)$ is equal to the maximal subspace of the form $\Circ(\mathbf{m})$.
	
	It remains to prove the claim. Fix $i \in \{ 0, \dots, n\}$ and $\mathbf{d} \in \Circ_i(V)$. Clearly the image of $\mathbf{d}$ under $C$ defines the same circulant space, hence $\Circ_i(V)$ is invariant. It remains to prove that it is an $i$-dimensional subspace.
	
	To achieve this, we review some combinatorics of polynomials. Let $Q = Q^{(n)} = \sum_{k = 0}^{n-1} q_kX^k \in \F_p[X]$ be a polynomial of degree $n-1$. In view of \cref{DerGGS:thm:KoenigRados}, we conduct Euclidean division by $(X-1)$, and write $Q^{(n)} = (X-1)Q^{(n-1)} + R_{n}(Q)$, for a polynomial $Q^{(n-1)}$ and a constant $R_{n}(Q) \in \F_p$. Iterating this, we write $Q^{(i)} = (X-1)Q^{(i-1)} + R_{i}(Q)$ for $i \in \{0, \dots, n-1\}$. Clearly $\deg Q^{(i)} = i-1$. The coefficients of $Q^{(i)}$ and the value of $R_i(Q)$ can be calculated in terms of the starting polynomial $Q^{(n)}$. We shall now perform this calculation. (Alternatively, we could compute the coefficients of $Q^{(n)}$ as a polynomial in $X-1$). Indeed, it is easy to check that
	\[
		Q^{(n-1)} = \sum_{k = 0}^{n-2} \sum_{\ell = k+1}^{n-1} q_{\ell} X^k \qund R_{n}(Q) = \sum_{k = 0}^{n-1} q_k.
	\]
	Thus, both the coefficient of $X^k$ in $Q^{(i)}$ and the value $R_i(Q)$ are weighted sums (i.e.\ positive $\F_p$-linear combinations) of the coefficients of $Q$. Write $\kappa(i,j,k)$ for the multiplicity of $q_j$ in the coefficient of $X^k$ in $Q^{(i)}$. The equation above shows that
	\[
		\kappa(i, j, k) = \sum_{\ell = k+1}^{n-1} \kappa(i+1, j, \ell).
	\]
	Consequently, for $i < n$, we find the familiar (at least when ignoring $j$) recursion formula
	\[
		\kappa(i, j, k) = \kappa(i, j, k+1) + \kappa(i+1, j, k+1),
	\]
	using that $\kappa(i,j,n) = 0$ for all $i,j \in \{0, \dots, n\}$. Since 
	$\kappa(n-1,j,k)$ is $1$ for $j > k$, and is equal to $0$ otherwise, we find
	\[
		\kappa(i, j, j + i - n) = 1 \qund \kappa(i, j, k) = 0 \quad\text{ for } k > j + i - n,
	\]
	and obtain the equality
	\[
		\kappa(i, j, k) = \binom{j-k-1}{n-i-1},
	\]
	where we agree on $\binom r s = 0$ for $r < s$. The remainder $R_i(Q)$ is equal to the sum of all coefficients of $Q^{(i)}$. Since $\deg Q^{(i)} = i-1$, we have to calculate the sum
	\[
		R_i(Q) = \sum_{k = 0}^{i-1} \sum_{j = 0}^{n-1} \kappa(i,j,k)q_j = \sum_{j = 0}^{n-1} q_j\sum_{k = 0}^{i-1} \binom{j-k-1}{n-i-1}.
	\]
	We may ignore all summands with $j < n-i$, since then $j-k-1 < n - i - 1$ for all $k \in \{0, \dots, i-1\}$, and the binomial coefficient is zero. Likewise we may ignore all cases where $k > j + i - n$. It remains to use `Stifel's law'. We find
	\[
		R_i(Q) = \sum_{j = n-i}^{n-1} q_j \sum_{k = 0}^{j+i-n} \binom{j-k-1}{n-i-1} = \sum_{j = n-i}^{n-1} q_j \sum_{k = n-i-1}^{j-1} \binom{k}{n-i-1} = \sum_{j = i}^n q_j \binom{j}{n-i}.
	\]
	Coming back to our circulant spaces, \cref{DerGGS:thm:KoenigRados} tells us that the condition $\mathbf{d} \in \Circ_i(V)$, for any $i$, translates to $R_j(E_{\mathrm{d}}) = 0$ for all $j \in \{i + 1, \dots, n\}$, and $R_{i}(E_{\mathrm{d}}) \neq 0$. By our computations, the map $R_j\colon V \to \F_p$ (for any $j \in \{1, n\}$) assigning to an element $\mathbf{d} = (d_0, \dots, d_{n-1})$, represented in the basis $\mathcal{B}$, the value of $R_j(E_{\mathrm{d}})$ as given above is $\F_p$-linear. Define $R \colon V \to \F_p^{p}$ by $\mathbf{d}R = (\mathbf{d}R_1, \dots, \mathbf{d}R_{n})$. This map is, due to $\mathcal B$ and the standard basis, represented by the matrix
	\[
		\begin{pmatrix}[1.5]
			\binom{0}{0} & \binom{1}{0} & \dots & \binom{n-2}{0} & \binom{n-1}{0}\\
			 & \binom{1}{1} & \dots & \binom{n-2}{1} & \binom{n-1}{1}\\
			 &  & \ddots & & \vdots\\
			 &  0 &  & \binom{n-2}{n-2} & \binom{n-1}{n-2}\\
			 &  &  &  & \binom{n-1}{n-1} \\
		\end{pmatrix},
	\]
	a right-justified Pascal~triangle. The subspace $\Circ_i(V)$ is the kernel of the composition $R \circ \pi_{\leq i}$, where $\pi_{\leq n-i}$ denotes the projection to the first $n-i$ coordinates. Since $R$ has full rank, the image under this map has dimension $n-i$, whence the kernel has dimension $i$. Thus $\dim \Circ_i(V) = i$.
\end{proof}

It is not true that every defining tuple gives rise to a unique \GGS-group. In particular, multiples of a given $\mathbf{e}$ define the same group (as a subgroup of $\Aut(X^*)$). Furthermore, certain reorderings of the entries give isomorphic groups, which helps us to reduce the difficulty of our computations. We use the following characterisation.
\begin{theorem}\cite{Pet19}\label{DerGGS:thm:classification ggs}
	Let $G$ and $H$ be two \GGS-groups over the $p$-regular tree defined by $\mathbf{e}$ and $\mathbf{d}$, respectively. Then the following two statements are equivalent:
	\begin{enumerate}
		\item $G \cong H$;
		\item there exist $\lambda, \mu \in \F_p^\times$ such that $\mathbf{e}_i = \mu \cdot \mathbf{d}_{\lambda \cdot i}$ for all $i \in \{1, \dots, p-1\}$.
	\end{enumerate}
\end{theorem}
This allows us to choose defining tuples with desirable properties.
\begin{corollary}\label{DerGGS:cor:standard form for e}
	Let $G$ be a \GGS-group. Then
	\begin{enumerate}
		\item there is an \GGS-group $G_{\mathbf{e}}$ isomorphic to $G$ such that $\mathbf{e}_1 = 1$, and
		\item there is an \GGS-group $G_{\mathbf{e}}$ isomorphic to $G$ such that $\mathbf{e'}_i = 1$ for some $i \in \{1, \dots, p-1\}$.
	\end{enumerate}
\end{corollary}

\subsection{Properties and structure of \GGS-groups}

We shall fix some further notation. Recall that, working with a given \GGS-group $G$, we shall denote the rooted generator $a$, the directed generator $b$, and we write $c$ for the commutator $[b, a]$. Furthermore we shall use the following shorthand notation for the conjugates of $c$,
\[
	c_i = c^{a^{i}} = [b^{a^{i}}, a].
\]
In particular, $c_0 = c$. We now describe the sections of the elements $c_i$. We will use this computation often and without constant reference. The sections of $c_i$ are the sections of $c_0$ cyclically shifted; in general, for any $g \in \Aut(X^*)$ and for any $i \in \Z$, the first level sections of $g^{a^{i}}$ are the sections of $g$, permuted by the inverse $a^{-i}$, since the sections of $a$ are trivial and
\[
	g^{a^{i}}|_j = a^{-i}|_j g|_{j^{a^{i}}} a^{i}|_{j^{a^{i}}} = g|_{j-i},
\]
i.e.\ $\psi(g^{a^{i}}) = \psi(g)^{a^{-i}}$, with $a^{-i}$ acting as a permutation of the index set $X$.

It is well-known that all \GGS-groups posses strong `self-similarity' properties, which is one reason making this class of groups an interesting object to study. We collect some statements into a lemma, for proofs see e.g.\ \cite{Pet19}.

\begin{lemma}[Fractality properties of \GGS-groups]\label{DerGGS:lem:fractality properties}
	Let $G$ be a \GGS-group. Then $G$ is self-similar, fractal and contracting with nucleus $\langle a \rangle \cup \langle b \rangle$.
\end{lemma}

Furthermore, due to Fern{\'a}ndez-Alcober and Zugadi-Reizabal, every \GGS-group is a regular branch group (see below). The same is not true for the analogues of \GGS-groups defined by constant tuples, explaining their divergent behaviour.

\begin{theorem}[Branching properties of \GGS-groups, cf.\ \cite{FZ13}]\label{DerGGS:thm:branching properties}
	Let $G$ be a \GGS-group with defining tuple $\mathbf{e}$. Then
	\begin{enumerate}
		\item $\psi(\gamma_3(\Stab_G(1))) = \gamma_3(G) \dottimes{p} \gamma_3(G)$, and
		\item $[G' \dottimes{p} G' : \psi(\Stab_G(1)')] = p^{\sym(\mathbf e)}$.
	\end{enumerate}
	Also $\Stab_G(2) \leq \gamma_3(G)$. In particular, $G$ is regular branch over $\gamma_3(G)$, and it is regular branch over $G'$ if $\mathbf{e}$ is non-symmetric.
\end{theorem}

This allows the application of the following lemma of Šuniḱ, which provides the analogue of \cref{DerGGS:thm:higher derived subgps} for level stabilisers.

\begin{lemma}\label{DerGGS:lem:stabs decompose}\cite[Lemma~10]{Sun07}
	Let $G$ be a regular branch group over $K \leq G$, such that $\Stab_G(n) \leq K$ for some $n \in \N$. Then for all $m \geq n$
	\[
		\psi(\Stab_G(m+1)) = \Stab_G(m) \dottimes{p} \Stab_G(m).
	\]
\end{lemma}

For \GGS-groups in particular, the value of $n$ is $2$.

We now begin with a study of certain small quotients of \GGS-groups which will play a role in the determination of the derived series. Most of these results are known, but using the next lemma, we give new and short proofs that have the benefit of being easily generalised to larger families of groups, as they (for the most part) do not involve $p$-group methods.

\begin{lemma}[Spine-counting lemma]\label{DerGGS:lem:b expsum is homomorphism}
	Let $G$ be a \GGS-group, let $g \in G$, and let $w(a,b)$ be a word in the letters $a$ and $b$ that evaluates to $g$ in $G$. Then clearly $w(1, b)$ evaluates to an element in $\langle b \rangle \cong \F_p$. The map
	\[
		\varepsilon: G \to \F_p, \quad g \to w(1, b)
	\]
	that assigns to $g$ its \emph{$b$-exponent sum} is well-defined and a homomorphism. In other words, the quotient of $G$ by the normal closure of $\langle a \rangle$ is a group of order $p$.
\end{lemma}

\begin{proof}
	Rewrite $w(a,b)$ into a word of the form
	\[
		{b^{a^{k_0}}}^{\ell_0} {b^{a^{k_1}}}^{\ell_1} \dots {b^{a^{k_{n-1}}}}^{\ell_{n-1}} a^{k_n},
	\]
	with \(\ell_i, k_i \in \F_p\) for $i \in \{0, \dots, n\}$ and \(k_i \neq k_{i+1}\) for \(i \in \{0, \dots, n - 2\}\). We may also assume \(n \geq 1\), since the statement is clearly true for $n = 0$. Let \(x\) be any letter of \(X\). Since all conjugates of \(b\) fix \(x\), a word representing \(w(a,b)|_x\) is
	\[
		b|_{x-k_0}^{\ell_0} b|_{x-k_1}^{\ell_1} \dots b|_{x-k_{n-1}}^{\ell_{n-1}}.
	\]
	Look at a length-$2$ subword \(b|_{x-k_j}^{\ell_j} b|_{x-k_{j+1}}^{\ell_{j+1}}\). Since 
	\(k_j \not = k_{j+1}\), at least one of the two sections is a power of $a$, since only one section of \(b\) is not. We collect all resulting elements of \(\langle a \rangle\) and combine consecutive powers of $b$ (which does not change the $b$-exponent sum) and obtain a word $w_x(a,b)$ with fewer syllables representing \(w(a,b)|_x\), such that $\sum_{x = 0}^{p-1} w_x(1,b) = w(1, b)$. Since $G$ is contracting (which is an consequence of the same argument), for every $w(a,b)$ there is a level $X^n$ for $n \in \N$, such that all sections of $w(a,b)$ at vertices $v \in X^n$ are in the nucleus of $G$, hence powers of $a$ or $b$. Clearly, all sections are minimal words, and the sum of their $b$-exponent sums is the $b$-exponent sum of $w(a,b)$. Thus this value does not depend on the word representing a certain element. Consequently $\varepsilon$ is well-defined and a homomorphism.
\end{proof}

Using this lemma, we are able to quickly recover some well-known facts about certain small quotients of \GGS-groups.

\begin{lemma}\label{DerGGS:lem:small quotients}
	Let $G$ be a \GGS-group. Then
	\begin{enumerate}
		\item $G/G'$ is isomorphic to an elementary abelian $p$-group of rank $2$,
		\item $G/\gamma_3(G)$ is isomorphic to the Heisenberg group over $\F_p$,
		\item $G/\Stab(1)'$ has order $p^{p+1}$, and
		\item the subgroup generated by the $p$\textsuperscript{th} powers of $G$ is contained in $G'$.
	\end{enumerate}
\end{lemma}

\begin{proof}
	(i): The quotient $G/G'$ is generated by the images of $a$ and $b$, hence a quotient of an elementary abelian $p$-group of rank $2$. Clearly $b$ is not a power of $a$, and $a$ is not a stabiliser of the first level, hence not in the commutator subgroup. Since the $b$-exponent sum of any element in the commutator subgroup is $0 \bmod p$, by \cref{DerGGS:lem:b expsum is homomorphism} the generator $b$ is not in $G'$ either.
	
	(ii): The centre of $G/\gamma_3(G)$ is generated by $c = [b,a]$. Since $c^p \equiv_{\gamma_3(G)} [b,a^{p}] = 1$, the group $G/\gamma_3(G)$ is isomorphic to the Heisenberg group over $\F_p$.
	
	(iii): The the first level stabiliser is generated by the conjugates of $b$ by powers of $a$. Since the position of the section equal to $b$ is different for every different power of $a$, it is generated by $p$ elements of order $p$, all of $b$-exponent sum $1$. Thus, as in (i), the quotient $\Stab_G(1)/\Stab_G(1)'$ is an elementary abelian $p$-group of rank $p$. Since $G/\Stab_G(1)$ is cyclic of order $p$, the result follows.
	
	(vi): Lastly, it is enough to prove that $b$ is not a $p$\textsuperscript{th} power. But the $b$-exponent sum of a $p$\textsuperscript{th} power is $0 \bmod p$, hence $b$ is not a $p$\textsuperscript{th} power by \cref{DerGGS:lem:b expsum is homomorphism}.
\end{proof}

\begin{proposition}\label{DerGGS:prop:gamma3 stab in G''}
	Let $G$ be a \GGS-group. Then
	\[
		[\Stab_G(1)', G'] = \gamma_3(\Stab_G(1)).
	\]
	In particular, $\gamma_3(\Stab_G(1)) \leq G''$.
\end{proposition}

\begin{proof}
	Recall from \cref{DerGGS:thm:branching properties} that
	$$
		\psi(\gamma_3(\Stab_G(1))) = \gamma_3(G) \dottimes{p} \gamma_3(G).
	$$
	The inclusion $[\Stab_G(1)', G'] \leq \gamma_3(\Stab_G(1))$ holds vacuously. We have to establish the other inclusion. Using \cref{DerGGS:prop:normal generation}, it is enough to prove that $(0: [c, b],\; \gnrc: \id)$ and $(0: [c, a],\; \gnrc: \id) \in \psi([\Stab_G(1)', G'])$, since $\gamma_3(G)$ is normally generated by $[c, b]$ and $[c, a]$. We distinguish two cases.\\
	
	\noindent\emph{Case 1}: The defining tuple is non-symmetric. Then by \cref{DerGGS:thm:branching properties}(ii)
	$$
		\psi(\Stab_G(1)') = G' \dottimes{p} G',
	$$
	and $(0: c,\; \gnrc: \id) \in \psi(\Stab_G(1)')$. By \cref{DerGGS:cor:standard form for e}(ii) we may assume $\mathbf{e'}_i = 1$ for some $i \in \{2, \dots, p-1\}$. Consequently $c|_i = a$, hence
	\begin{align*}
		c_{p-i}|_0 &= c|_i = a,\\
		(c_{p-i}^{\mathbf{e}_{p-1}}c^{-1})|_0 &= a^{\mathbf{e}_{p-1}} a^{-\mathbf{e}_{p-1}} b = b.
	\end{align*}
	Since $c_{p-i}$ and $c$ are elements of $G'$, we obtain that
	\begin{align*}
		[(0: c,\; \gnrc:\id), \psi(c_{p-i})] &= (0: [c, a],\; \gnrc:\id), \quad\text{ and}\\
		[(0: c,\; \gnrc:\id), \psi(c_{p-i}^{\mathbf{e}_{p-1}}c^{-1})] &= (0:[c, b],\; \gnrc:\id)
	\end{align*}
	are both elements of $\psi([\Stab_G(1)', G'])$.\\
		
	\noindent\emph{Case 2}: The defining tuple is symmetric. Since it is not constant, the prime $p$ is necessarily greater then $3$. By \cref{DerGGS:cor:standard form for e}(i), we may assume $\mathbf{e}_1 = \mathbf{e}_{p-1} = 1$. Observe that
	\[
		\psi([b, b^a]) = [\psi(b), \psi(b^a)] = (0: [b, a],\; 1: [a, b],\; \gnrc: \id) = (0: c,\; 1: c^{-1},\; \gnrc: \id).
	\]
	Now let $j \in \F_p^\times\smallsetminus\{1, p-1\}$, which exists since $p > 3$. Then
	\[
		\psi([a^{2j}, b]) = (0: a^{-\mathbf{e}_{p-2j}}b,\; 2j: b^{-1}a^{\mathbf{e}_{2j}},\; \gnrc: a^{\mathbf{e}_{\gnrc} - \mathbf{e}_{\gnrc-2j}}),
	\]
	in particular $[a^{2j},b]|_{j} = a^{\mathbf{e}_{j} - \mathbf{e}_{p-j}} = \id$, since $\mathbf{e}$ is symmetric. Thus
	\begin{align*}
		[a^{2j}, b]^{a^{1-j}}|_1 &= [a^{2j}, b]|_j = \id, \quad\text{ and }\\
		[a^{2j}, b]^{a^{1-j}}|_0 &= [a^{2j}, b]|_{j-1} = a^{\mathbf{e}_{j-1}-\mathbf{e}_{p-j-1}},
	\end{align*}
	since $j \not\in \{1, p-1\}$. If there is such a $j$ with $\mathbf{e}_{j-1} \neq \mathbf{e}_{p-j-1}$, let $i \in \Z$ be such that $i \equiv_p (\mathbf{e}_{j-1}-\mathbf{e}_{p-j-1})^{-1}$, and observe that
	\begin{align*}
		\psi([[b^{a}, b], [a^{2j}, b]^i]) &= [(0: c,\; 1: c^{-1},\; \gnrc: \id), (0: a,\; 1: \id,\; \gnrc: [a^{2j}, b]^i)|_{\gnrc}]\\
		&= (0:[c,a],\; \gnrc:\id).
	\end{align*}
	But such an element $j$ must always exist. Assume otherwise, for contradiction. Then for any $j \notin \{1, p-1\}$
	\begin{align*}
		\mathbf{e}_j = \mathbf{e}_{(1 + j) - 1} = \mathbf{e}_{p - (1 + j) - 1} = \mathbf{e}_{j + 2}.
	\end{align*}
	But since every element of $\F_p^\times$ is a multiple of $2$, this implies that $\mathbf{e}$ is constant.
	
	It remains to show $(0:[c, b],\; \gnrc: \id) \in \psi([\Stab_G(1)', G'])$. If $\mathbf{e}_{p-2} \neq 0$, consider the element $g = [a^{2j}, b]^{i \mathbf{e}_{p-2}}[a^{2}, b]$, which fulfils
	\begin{align*}
		g|_0 = a^{\mathbf{e}_{p-2}} a^{-\mathbf{e}_{p-2}} b \quad\text{ and }\quad
		g|_1 = a^{\mathbf{e}_1 - \mathbf{e}_{p-1}} = \id,
	\end{align*}
	using the fact that $\mathbf{e}$ is assumed to be symmetric.
	If otherwise $\mathbf{e}_{p-2} = 0$, set $g = [a^{2}, b]$, which has the sections
	\begin{align*}
		g|_0 = a^{-\mathbf{e}_{p-2}}b = b, \quad\text{ and }\quad g|_1 = a^{\mathbf{e}_1 - \mathbf{e}_{p-1}} = \id.
	\end{align*}
	In both cases
	\[
		\psi([[b^{a}, b], g]) = (0:[c,b],\; \gnrc:\id). \qedhere
	\]
\end{proof}

Using \cref{DerGGS:prop:gamma3 stab in G''}, we may derive a nice corollary.

\begin{corollary}\label{DerGGS:cor:branch over G''}
	Let $G$ be a \GGS-group. Then $G$ is branch over $G''$, independent on the value of $\sym(\mathbf{e})$.
\end{corollary}

\begin{proof}\belowdisplayskip=-11pt
	Using \cref{DerGGS:thm:branching properties}(i) and \cref{DerGGS:prop:gamma3 stab in G''} we find the inclusion
	\begin{align*}
		G'' \dottimes{p} G'' \leq \gamma_3(G) \dottimes{p} \gamma_3(G) &= \psi(\gamma_3(\Stab(1))) \\&= \psi([\Stab(1)', G']) \leq \psi(G'').
	\end{align*}\qedhere
\end{proof}

\section{The derived series of \GGS-groups}
The main difficulty for the computation of the index $|G:G^{(n)}|$ for all $n > 1$ is the calculation of $|G:G''|$. We now begin with this case, applying the theory of circulant spaces sketched in \cref{DerGGS:prop:flag of circulant spaces}. Recall that we can compute the circulant space generated by some vector $\mathbf{d} \in \F_p^p$ by counting the number of trailing zeros in the image of $\mathbf{d}$ under the linear map $R$, and that the last three components may be computed as
\[
	\mathbf{d}R_p = \sum_{i = 0}^{p-1} d_i, \quad \mathbf{d}R_{p-1} = \sum_{i = 0}^{p-1} id_i,\qund \mathbf{d}R_{p-2} = \sum_{i = 0}^{p-1} \binom{i}{2} d_i.
\]

\begin{proposition}\label{DerGGS:prop:G'' index}
	Let $G$ be a \GGS-group with defining tuple $\mathbf{e}$. Then
	\[
		\log_p |G:G''| = p + 1 + \con(\mathbf e') + \sym(\mathbf e'') - \sym(\mathbf e).
	\]
\end{proposition}

\begin{proof}
	Since $G' \leq \Stab_G(1)$, the second derived subgroup $G''$ is contained in $\Stab_G(1)'$, which is in turn contained in $G' \dottimes{p} G'$ by \cref{DerGGS:thm:branching properties}. By \cref{DerGGS:prop:gamma3 stab in G''} and \cref{DerGGS:thm:branching properties} we find
	\[
		\psi(G'') \geq \psi(\gamma_3(\Stab(1))) = \gamma_3(G) \dottimes{p} \gamma_3(G),
	\]
	so that $G''$ is wedged in between the subgroups $\psi^{-1}(\gamma_3(G) \dottimes{p} \gamma_3(G))$ and $\psi^{-1}(G' \dottimes{p} G')$. Since the quotient $G'/\gamma_3(G)$ is cyclic of order $p$ by \cref{DerGGS:lem:small quotients}(i) and (ii), we may identify the quotient
	\[
		(G' \dottimes{p} G')/(\gamma_3(G) \dottimes{p} \gamma_3(G)) \cong G'/\gamma_3(G) \dottimes{p} G'/\gamma_3(G)
	\]
	with a $p$-dimensional $\F_p$-vector space $V$, such that $G'' \cdot (\gamma_3(G) \dottimes{p} \gamma_3(G))$ represents a sub-vector space $W$ of $V$. We calculate the dimension of this subspace. Knowing it, we can easily deduce the index of the second derived subgroup. Recall that $\log_p [G' \dottimes{p} G' : \Stab(1)'] = \sym(\mathbf{e})$ by \cref{DerGGS:thm:branching properties}, and has that $\Stab(1)'$ has $p$-logarithmic index ${p+1}$ in $G$ by \cref{DerGGS:lem:small quotients}. Thus
	\begin{align*}
		\log_p |G:G''| &= \log_p |G:\Stab_G(1)'| + \log_p [\Stab_G(1)': G'']\\
		&= \log_p |G:\Stab_G(1)'| + \log_p |G' \dottimes{p} G': G''|\\
		&\quad\quad\quad - \log_p |G' \dottimes{p} G' : \Stab_G(1)'|\\
		&= p + 1 + \codim_V(W) - \sym(\mathbf{e}).
	\end{align*}
	Note that it is a consequence of \cite[Lemma~3.5]{FZ13} (and indeed, of the calculation of $[b, b^a]$ at the beginning of `Case 2' in the proof of \cref{DerGGS:prop:gamma3 stab in G''}) that the subspace $U$ represented by $\Stab(1)'$ in $V$ is equal to $\Circ_{p-1}(V)$ if $\mathbf{e}$ is symmetric. We will not use this fact here, but find it helpful to keep in mind.

	To compute the dimension of $W$, we calculate the images of a set of generators of $W$ in $V$. Since $G'$ is generated by the elements $\{c^g \mid g \in G\}$, the second derived subgroup $G''$ is generated by elements of the form
	\[
		[c^{g_1}, c^{g_2}]^{g_3} = [c, c^{g_2g_1^{-1}}]^{g_1g_3}
	\]
	for some $g_1, g_2$ and $g_3 \in G$, i.e.\ a generating set is given by
	\[
		\{ [c, c^g]^h \mid g, h \in G\}.
	\]
	Since we calculate modulo $\gamma_3(G)\dottimes{p}\gamma_3(G)$, we may restrict the choices of $g$ and $h$ significantly. Indeed, writing $g = a^i s$ and $h = a^j t$ for some $i, j \in F_p$ and $s, t \in \Stab(1)$, we see that, for any $x \in X$ (and writing $y = x^{a^{j}}$ and $z = y^{a^i} = x^{a^{i+j}}$), and using the fact that the section map is a homomorphism on $\Stab(1)$,
	\[
		[c, c^{a^is}]^{a^jt}|_x = [c, c^{a^is}]^{t}|_{y} = [c|_{y}, c^{a^is}|_{y}]^{t|_{y}} = [c|_{y}, (c|_z)^{s|_z}]^{t|_{y}} \equiv_{\gamma_3(G)} [c|_{y}, c|_z].
	\]
	Thus $s$ and $t$ play no role modulo $\gamma_3(G) \dottimes{p} \gamma_3(G)$, and we may restrict to the generating set
	\[
		\{ [c, c^{a^{j-i}}]^{a^{-j}} \mid i, j \in \F_p \} = \{ [c_i, c_j] \mid i, j \in \F_p \}.
	\]
	We do also know how $a$ acts (via the conjugation action on $\Stab(1)'$) on the vector space $U$. To describe this action, we choose a basis; since $G'/\gamma_3(G)$ is generated by the image $\overline{c}$ of $c$, we decide to consider this abelian quotient as a group written additively, i.e. $c^i \gamma_3(G) = i\overline{c}$ for any $i \in \F_p$. Our basis is hence the standard basis with respect to $\overline{c}$
	\[
		\left( (\overline{c}, 0, \dots, 0), (0, \overline{c}, 0, \dots, 0), \dots, (0, \dots, 0, \overline{c}) \right)
	\]
	Now conjugation by $a$ (on any element $s$ of the stabiliser) cyclically permutes the entries of the tuple $\psi_1(s)$. This corresponds to the linear map cyclically permuting the basis elements, hence extends naturally to an action on the full space $V$. Since $W$ is the subspace representing $G''$, a normal subgroup, it is invariant under this action, i.e.\ it is a circulant space. Notice that $[c_i, c_j]^{a^{-k}} = [c_{i+k}, c_{j+k}]$. Furthermore
	\[
		[c, c_i]^{-1} = [c_i, c] = [c^{a^i}, c] = [c, c_{p-i}]^{a^i},
	\]
	thus the set $\{[c, c_i] \mid i \in \{ 1, \dots, {(p-1)}/2 \} \}$ normally generates $G''$. Equivalently, the images under $\theta$ (the map $\Stab(1) \ni g \mapsto (\log_a(g|_0), \dots, \log_a(g|_{p-1}))$ defined at the beginning of \cref{sec:def tuples})
	\[
		\mathbf{d_i} = \theta([c, c_i])
	\]
	generate $W$ as a cyclically invariant subspace, i.e.\
	\[
		W = \Circ\left(\,\left\{ \mathbf{d_i} \;\middle|\; i \in \left\{ 1, \dots, {(p-1)}/2 \right\}\, \right\}\,\right).
	\]
	By \cref{DerGGS:prop:flag of circulant spaces}, it is enough to compute the rank of $\Circ(\mathbf{d_i})$ for every $i$. So let us compute the images of $[c, c_i]$ under $\theta$. We begin with the cases $i \in \{ 2, \dots, {(p-1)}/2 \}$, where we find
	\begin{align*}
		\psi([c, c_i]) &=  \left(\begin{array}{rlrl}
			0:& [b^{-1}a^{\mathbf{e}_{p-1}}, a^{\mathbf{e'}_{p-i}}], & 1:& [a^{-\mathbf{e}_{1}}b, a^{\mathbf{e'}_{p-i+1}}]\\
			i:& [a^{\mathbf{e'}_{i}}, b^{-1}a^{\mathbf{e}_{p-1}}]), & i+1:& [a^{\mathbf{e'}_{i+1}}, a^{-\mathbf{e}_{1}}b]\\
			\gnrc:& \id &&
		\end{array}\right)\\
		&\equiv
		\left(\begin{array}{rlrl}
			0 :& [b, a]^{-\mathbf{e'}_{p-i}}, & 1 :& [b, a]^{\mathbf{e'}_{p-i+1}},\\
			i :& [b, a]^{\mathbf{e'}_{i}}, & i+1 :& [b, a]^{-\mathbf{e'}_{i+1}},\\
			\gnrc :& \id &&
		\end{array}\right)\pmod{{\gamma_3(G) \dottimes{p} \gamma_3(G)}}.
	\end{align*}
	By our choice of basis for $V$, we find
	\begin{equation*}\label{eq:yi}
		\mathbf{d_i} = \theta([c, c_i]) = (0: -\mathbf{e'}_{p-i},\; 1: \mathbf{e'}_{p-i+1},\; i: \mathbf{e'}_{i},\; i+1: -\mathbf{e'}_{i+1}).\tag{$\ast_i$}
	\end{equation*}
	The last generator $[c, c_1]$ has a slightly different form,
	\begin{align*}
		\psi([c, c_1]) &= (0: [b^{-1}a^{\mathbf{e}_{p-1}}, a^{\mathbf{e'}_{p-1}}]),\; 1: [a^{-\mathbf{e}_{1}}b, b^{-1}a^{\mathbf{e}_{p-1}}], \; 2: [a^{\mathbf{e'}_{2}}, a^{-\mathbf{e}_{1}}b],\; \gnrc: \id)\\
		&\equiv(0: c^{-\mathbf{e'}_{p-1}},\; 1: c^{\mathbf{e}_{p-1}-\mathbf{e}_{1}},\; 2: c^{-\mathbf{e'}_2},\; \gnrc : \id) \pmod{{\gamma_3(G) \dottimes{p} \gamma_3(G)}},
	\end{align*}
	and gives rise to the vector
	\begin{equation*}\label{eq:y1}
		\mathbf{d_1} = \theta([c, c_1]) = (0: -\mathbf{e'}_{p-1},\; 1: \mathbf{e}_{p-1} - \mathbf{e}_1,\; 2: -\mathbf{e'}_{2},\; \gnrc: 0).\tag{$\ast_1$}
	\end{equation*}
	Using these descriptions, we are able to compute some of the circulant spaces associated to $\mathbf{d_i}$ by applying the maps $R_p, R_{p-1}$ and $R_{p-2}$ to $\mathbf{d_i}$. To do so, we make a case distinction.\\
	
	\noindent\emph{Case 1}: If the second differences vector $\mathbf{e''}$ is not symmetric, there is some index $i \in \{3, \dots, p-1\}$ such that $\mathbf{e''}_{i+1} \neq \mathbf{e''}_{p - i + 1}$. In this case, the sum of the components of the corresponding vector $\mathbf{d_i}$ is
	$$
		\mathbf{d_i}R_p = (\mathbf{e'}_{p-i} - \mathbf{e'}_{p-i+1}) - (\mathbf{e'}_{i} - \mathbf{e'}_{i+1}) = \mathbf{e''}_{p-i+1} - \mathbf{e''}_{i+1} \neq 0,
	$$
	and by \cref{DerGGS:thm:KoenigRados} the rank of the circulant matrix $\Circ(\mathbf{d_i}) = p$, since the sum of the coefficients of $\mathbf{d_i}$ is not $0$. The subspace $W$ is equal to the full space $V$, and we find
	\[
		\log_p |G:G''| = \log_p |G: G' \dottimes{p} G'|.
	\]
	Assume that the first differences vector $\mathbf{e'}$ is constant. Then $\mathbf{e''}$ is constant $0$, hence symmetric, a contradiction. Also, by \cref{DerGGS:lem:e sym implies snd der sym}, the asymmetricality of $\mathbf{e''}$ forces $\mathbf{e}$ to not be symmetric, hence $G' \dottimes{p} G'$ is in fact equal to $\Stab(1)'$ and a subgroup of $G$. We conclude that the index of the second derived subgroup fulfils
	\[
		\log_p |G : G''| = p + 1 = p + 1 + \underbrace{\con(\mathbf e')}_{=0} + \underbrace{\sym(\mathbf e'')}_{=0} - \underbrace{\sym(\mathbf e)}_{=0}.
	\]
	
	\noindent\emph{Case 2}: Now assume that the defining tuple $\mathbf{e''}$ is symmetric. Using the description of the generators given in the first case, we see that $\Circ(\theta([c, c_i])) \subseteq \Circ_{p-1}(V)$. We furthermore need to consider the generator $[c, c_1]$, that played no role in the previous case. Since
	\begin{align*}
		\mathbf{d_1} R_p = -\mathbf{e'}_{p-1} + \mathbf{e}_{p-1} - \mathbf{e}_{1} - \mathbf{e'}_2 = 2(\mathbf{e}_{p-1} - \mathbf{e}_{1}) + (\mathbf{e}_{2} - \mathbf{e}_{p-2}) = 0
	\end{align*}
	by \cref{DerGGS:lem:linear eq a consequence of snd der sym}, the sum of the coefficients of $[c, c_1]$ vanishes, and $\Circ([c, c_1])$ is contained in $\Circ_{p-1}(V)$ as well. Consequently, $W$ is a subspace of $\Circ_{p-1}(V)$. We can go on and reduce the possible size of $\Circ([c, c_1])$ further. Consider the image of the cyclic shift $[c, c_1]^{a^{-1}}$ under $\theta$,
	\[
		\widehat{\mathbf{d_1}} = \theta([c, c_1]^{a^{-1}}) = (\mathbf{e}_{p-1}-\mathbf{e}_1, -\mathbf{e'}_1, 0, \dots, 0, -\mathbf{e'}_{p-2}).
	\]
	Of course, this vector generates the same circulant space. But clearly
	\[
		\widehat{\mathbf{d_1}} R_1 = 0\cdot(\mathbf{e}_{p-1}-\mathbf{e}_1) - \mathbf{e'}_1 - (p-1) \mathbf{e'}_{p-2} = - \mathbf{e'}_1 + \mathbf{e'}_{p-2} = 0,
	\]
	since $\mathbf{e'}$ is symmetric. Hence $\Circ(\mathbf{d_1}) \subseteq \Circ_{p-2}(V)$.
	
	Contrary to the situation in the previous case, the symmetry of $\mathbf{e''}$ does not force $\mathbf{e}$ or $\mathbf{e'}$ to be symmetric (resp.\ constant). Therefore, we have to make another case distinction.\\
	
	\noindent\emph{Subcase 2.1}: Assume additionally that $\mathbf{e'}$ is not constant, i.e.\ that $\mathbf{e''}$ is not zero. We want to show that $\Circ(\mathbf{d_i}) = \Circ_{p-2}(V)$ for some $i \in \{1, \dots, {(p-1)}/2\}$.	 
	 Notice that for all $i \in \{2, \dots, p-1\}$ we find
	\begin{align*}
		\mathbf{e'}_{i+1} - \mathbf{e'}_{p-i+1} = \mathbf{e''}_{i+1} + \mathbf{e'}_{i+2} - \mathbf{e'}_{p-i} - \mathbf{e''}_{p-i+1}
		= \mathbf{e'}_{i+2} - \mathbf{e'}_{p-i},
	\end{align*}
	since $\mathbf{e''}$ is symmetric, and consequently
	\begin{equation*}\label{eq:middle}
		\mathbf{e'}_{i+1} - \mathbf{e'}_{p-i+1}
		= \mathbf{e'}_{(p-1)/2} - \mathbf{e'}_{(p+1)/2} = \mathbf{e''}_{(p+1)/2}.\tag{$\dagger$}
	\end{equation*}
	Building on this calculation, we can compute
	\begin{align*}
		\mathbf{d_i} R_1 &= 
		\binom{0}{1} \cdot -\mathbf{e'}_{p-i} + \binom{1}{1}\mathbf{e'}_{p-i+1} + \binom{i}{1}\mathbf{e'}_{i} - \binom{i+1}{1}\mathbf{e'}_{i+1}\\
		&= \mathbf{e'}_{p-i+1} - \mathbf{e'}_{i+1} + i(\mathbf{e'}_{i} - \mathbf{e'}_{i+1})\\
		&= i \mathbf{e''}_{i+1} - \mathbf{e''}_{(p+1)/2}.
	\end{align*}
	This is not $0$ if $\mathbf{e''}_{(p+1)/2} = 0$, since there must exist some $i \in \{3, \dots, {(p-3)}/2\}$ such that $\mathbf{e''}_{i} \neq 0$, because $\mathbf{e''}$ is by assumption symmetric and non-zero. If otherwise $\mathbf{e''}_{(p+1)/2} \neq 0$, for $i = {(p-1)}/2$ the equation becomes
	\[
		\mathbf{d_i} R_{p-1} = i \mathbf{e''}_{i+1} - \mathbf{e''}_{(p+1)/2} = (p-3)/2 \cdot \mathbf{e''}_{(p+1)/2} \neq 0.
	\]
	Thus $\Circ_{p-1}(V) = \Circ(\mathbf{d_i}) \subseteq W \subseteq \Circ_{p-1}$, and we can verify that
	\begin{align*}
		\log_p |G : G''| &= p + 1 + \codim_V(W) - \sym(\mathbf{e})\\
		&= p + 2 - \sym(\mathbf{e})\\
		&= p + 1 + \underbrace{\con(\mathbf e')}_{=0} + \underbrace{\sym(\mathbf e'')}_{=1} - \sym(\mathbf e).
	\end{align*}
	
	\noindent\emph{Subcase 2.2}: Finally, assume, in addition to $\mathbf{e''}$ being symmetric, that $\mathbf{e'}$ is constant. Write $k \in \F_p^\times$ for the constant value of $\mathbf{e'}$. Note that $k \neq 0$, since otherwise $\mathbf{e'}$ is zero and $\mathbf{e}$ constant.
	
	As we have argued above, $W$ is contained in $\Circ_{p-1}(V)$. In view of the formula we want to establish, we aim to prove that $W = \Circ_{p-2}(V)$. We first show that $\Circ_{p-2}(V) \subseteq W$, by computing
	\begin{align*}
		\mathbf{d_1} R_{p-2} = -\binom{0}{2}\mathbf{e'}_{p-2} + \binom{1}{2}\left(\mathbf{e}_{p-1}-\mathbf{e}_{1}\right) - \binom{2}{2}\mathbf{e'}_1 = -\mathbf{e'}_1 = -k \neq 0,
	\end{align*}
	hence $\Circ_{p-2}(V) = \Circ(\mathbf{y}_1) \subseteq W$.
	
	It remains to show that $\Circ(\mathbf{d_i}) \subseteq \Circ_{p-2}(V)$ for all $i \in \{2, \dots, \frac{p-1}2\}$, i.e.\ to show that \(\mathbf{d_i} R_1 = 0\). We calculate
	\begin{align*}
		\mathbf{d_i} R_{p-1} &= 
		\binom{0}{1} \cdot -\mathbf{e'}_{p-i} + \binom{1}{1}\mathbf{e'}_{p-i+1} + \binom{i}{1}\mathbf{e'}_{i} - \binom{i+1}{1}\mathbf{e'}_{i+1}\\
		&= k + ik - (i+1)k = 0
	\end{align*}
	Thus $W = \Circ_{p-2}(V)$ and
	\begin{align*}
		\log_p |G : G''| = p + 3 - \sym(\mathbf{e}) = p + 1 + \underbrace{\con(\mathbf e')}_{=1} + \underbrace{\sym(\mathbf e'')}_{=1} - \sym(\mathbf e),
	\end{align*}
	which concludes the proof.
\end{proof}
In fact, our proof gives more than just the index, but allows for some structural description of $G''$. We record this in the following corollary.
\begin{corollary}\label{DerGGS:cor:structure of G''}
	Let $G$ be a \GGS-group. Write $V$ for the vector space $G'/\gamma_3(G) \dottimes{p} G'/\gamma_3(G)$. Write $\pi$ for the quotient map $G \to G/\gamma_3(\Stab(1))$. Then
	\[
		\psi(G'')
		= \left\{ (g_1, \dots, g_p) \in G' \dottimes{p} G' \mid (\pi(g_1), \dots, \pi(g_p)) \in \Circ_{p-i}(V) \right\},	
	\]
	where $i = \con(\mathbf{e}) + \sym(\mathbf{e''})-\sym(\mathbf{e})$.
\end{corollary}

To describe the series of iterated local laws, we record the following lemma.

\begin{lemma}\label{DerGGS:lem:pth powers}
	The group of $p$\textsuperscript{th} powers of $G'$ is contained in $\gamma_3(\Stab_G(1))$. In particular, it is contained in $G''$.
\end{lemma}

\begin{proof}
	The group $G'$ is geometrically contained in the direct product $G \dottimes{p} G$. Now $G/\gamma_3(G)$ is isomorphic to the Heisenberg group over $\F_p$, hence of exponent $p$. Consequently $(G')^p \leq \psi^{-1}(\gamma_3(G) \dottimes{p} \gamma_3(G)) = \gamma_3(\Stab_G(1))$.
\end{proof}

Now that we have established the value of $|G:G''|$, it remains derive our main results.

\begin{proof}[Proof of \cref{DerGGS:thm:higher derived subgps}]
	Assume that the given equation holds true for some $n \in \N$. Then we find
	\[
		\psi(G^{(n+1)}) = \psi(G^{(n)})' = (G^{(n-1)} \dottimes{p} G^{(n-1)})' = G^{(n)} \dottimes{p} G^{(n)},
	\]
	since $G^{(n)} \leq \Stab_G(1)$. Thus by induction it is enough to consider the case $n = 3$, or the case $n = 2$, respectively.
	
	First assume that $\con(\mathbf{e'}) + \sym(\mathbf{e''}) - \sym(\mathbf{e}) = 0$. By \cref{DerGGS:thm:branching properties} we find $\psi(\Stab_G(1)') = G' \dottimes{p} G'$, and by \cref{DerGGS:prop:G'' index} $\log_p[G: G''] = p + 1$. By \cref{DerGGS:lem:small quotients}(iii), we also have $\log_p |G:\Stab_G(1)'| = p + 1$, and since $G'' \leq \Stab_G(1)'$, the subgroups $G''$ and $\Stab(1)' = \psi^{-1}(G'\dottimes{p} G')$ coincide. Hence the equation holds for $n = 2$.
	
	Now we drop the assumption on the defining tuple, and prove the desired equation for $n = 3$. Since $G''$ is normally generated by the elements $\{[c, c^g]\mid g \in G\}$, we have to prove that $(0:[c, c^g],\; \gnrc: \id ) \in \psi(G^{(3)})$ for any $g \in G$, then an application of \cref{DerGGS:prop:normal generation} concludes the proof.
	
	Let $g \in G$ be arbitrary. Since $G$ is fractal, we may find $\widehat{g} \in \Stab(1)$ such that $\widehat{g}|_0 = g$. Furthermore, we know by \cref{DerGGS:cor:structure of G''} that we find
	\[
		h = (0: c,\; 1:c^{-2},\; 2:c,\; \gnrc: \id) \in \psi(G'').
	\]
	Thus
	\[
		[h, h^{\widehat g a^{p-2}}] = (\gnrc: [h|_{\gnrc}, h^{\widehat{g}}|_{\gnrc+2}]) = (0: [c, c^g],\; \gnrc: \id) \in \psi(G^{(3)}).\qedhere
	\]
\end{proof}

\begin{proof}[Proof of \cref{DerGGS:thm:main}]
	Using \cref{DerGGS:thm:higher derived subgps}, \cref{DerGGS:prop:G'' index} and \cref{DerGGS:thm:branching properties}, we obtain
	\begin{align*}
		\log_p |G'' : G^{(3)}| &= \log_p |G' \dottimes{p} G' : G'' \dottimes{p} G''| - \log_p |G' \dottimes{p} G': \psi(G'')|\\
		&= \log_p (|G' : G''|^p) - \con(\mathbf{e'}) + \sym(\mathbf{e''}) \\
		&= (p-1)(p + \con(\mathbf{e'}) + \sym(\mathbf{e''})) - p\sym(\mathbf{e})).
	\end{align*}
	Using \cref{DerGGS:thm:higher derived subgps} yet again, we find for $n > 2$
	\begin{align*}
		\log_p |G^{(n-1)}:G^{(n)}| &= \log_p |G^{(n-2)} \dottimes{p} G^{(n-2)} : G^{(n-1)} \dottimes{p} G^{(n-1)}|\\
		&= p \log_p |G^{(n-2)}:G^{(n-1)}|,
	\end{align*}
	and consequently
	\begin{align*}
		\log_p |G'' : G^{(n)}| &= \sum_{i = 2}^{n-1} \log_p |G^{(i)} : G^{(i+1)}|\\
		&= (p^{n-2}-1)/(p-1) \cdot \log_p |G'': G^{(3)}|\\
		&= (p^{n-2}-1)(p + \con(\mathbf{e'}) + \sym(\mathbf{e''})) - (p^{n-1}-p)/(p-1) \cdot \sym(\mathbf{e}),
	\end{align*}
	hence by \cref{DerGGS:prop:G'' index}
	\[
		\log_p |G : G^{(n)}| = p^{n-2}(p + \con(\mathbf{e'}) + \sym(\mathbf{e''})) - (p^{n-1}-1)/(p-1)\cdot\sym(\mathbf{e}) + 1.\qedhere
	\]
\end{proof}

\begin{corollary}\label{DerGGS:cor:Derived and iterated local}
	Let $G$ be a \GGS-group defined by a non-constant vector. Then the series of iterated local laws $(\mathrm L_n(G))_{n \in \N}$ coincides with the derived series $(G^{(n)})_{n \in \N}$.
\end{corollary}

\begin{proof}\belowdisplayskip=-11pt
	Let $G$ be a \GGS-group. We have to show $G^{(n)} = \mathrm L_n(G)$ for all $n \in \N$. It is enough to show that all $p$\textsuperscript{th} powers in $G^{(n)}$ are contained in $G^{(n+1)}$ for all $n \in \N$; for in this case all words involving a $p$\textsuperscript{th} power give rise to verbal subgroups that are contained in the $(n+1)$\textsuperscript{st} derived subgroup.
	
	By \cref{DerGGS:lem:small quotients}(iv) and \cref{DerGGS:lem:pth powers} this holds for $n \in \{1, 2\}$. Hence we assume that $n > 2$. Then, by \cref{DerGGS:thm:higher derived subgps}, the group $\psi(G^{(n)})$ is contained in $G^{(n-1)} \dottimes{p} G^{(n-1)}$. Using induction, we see that	
	\begin{align*}
		\psi(G^{(n)})^p &= (G^{(n-1)} \dottimes{p} G^{(n-1)})^p\\
		&\leq ((G^{(n-1)})^p \dottimes{p} (G^{(n-1)})^p)\\
		&\leq (G^{(n)} \dottimes{p} G^{(n)}).
	\end{align*}\qedhere
\end{proof}


\providecommand{\bysame}{\leavevmode\hbox to3em{\hrulefill}\thinspace}
\providecommand{\MR}{\relax\ifhmode\unskip\space\fi MR }
\providecommand{\MRhref}[2]{%
  \href{http://www.ams.org/mathscinet-getitem?mr=#1}{#2}
}
\providecommand{\href}[2]{#2}

\end{document}